\documentclass{amsart}
\usepackage{url}
\usepackage{graphicx, float}

\usepackage{mathtools}
\usepackage{hyperref}

\newtheorem{theorem}{Theorem}[section]

\theoremstyle{definition}

\theoremstyle{remark}

\numberwithin{equation}{section}

\begin{document}
\title[AMS Article Template]{A Proof of The Triangular Ashbaugh--Benguria--Payne--P\'{o}lya--Weinberger Inequality}


\author{Ryan Arbon}
\author{Mohammed Mannan}
\author{Michael Psenka}
\author{Seyoon Ragavan}

\address{Department of Mathematics, Princeton University, Princeton, NJ 08544, USA}
\email{rarbon@princeton.edu, mmannan@princeton.edu, mpsenka@princeton.edu, sragavan@princeton.edu}
\curraddr{}

\begin{abstract}
In this paper, we show that for all triangles in the plane, the equilateral triangle maximizes the ratio of the first two Dirichlet--Laplacian eigenvalues. This is an extension of work by Siudeja \cite{Siudeja}, who proved the inequality in the case of acute triangles. The proof utilizes inequalities due to Siudeja and Freitas \cite{quad}, together with improved variational bounds.
\end{abstract}

\maketitle

\pagestyle{plain}

\bibliographystyle{amsplain}

\section{Introduction}

For triangles in the Euclidean plane, the explicit values for the eigenvalues of the Dirichlet--Laplacian problem are only known in the case of the equilateral, 30-60-90, and 45-45-90 triangles. However, it is known that for a given domain $D$ in the plane, the Dirichlet--Laplacian eigenvalues form a non-decreasing sequence, which we order as $\{\lambda_i\}_{i \in \mathbb{N}}$. From now on, given a domain $D$ in the plane, we will use the phrase ``the eigenvalues of $D$'' to refer to the Dirichlet--Laplacian eignvalues of $D$. 

The Payne--P\'{o}lya--Weinberger (PPW) inequality dates back to 1955, when L.\ Payne, G.\ P\'{o}lya, and H.\ Weinberger published a paper \cite{PPW} proving a bound on the ratio of the first two eigenvalues $\lambda_2 / \lambda_1$ of a bounded domain $D$ in the plane, namely that $\lambda_2 / \lambda_1 \le 3$. Payne, P\'{o}lya, and Weinberger conjectured that this ratio is maximized when $D$ is the disc, that is:

\begin{equation}
    \frac{\lambda_2}{\lambda_1} \leq \left. \frac{\lambda_2}{\lambda_1} \right\vert_{\text{disc}} \approx 2.539.
\end{equation}

The original PPW inequality was generalized to dimension $n$ by Thompson in \cite{thompson}, who showed that

\begin{equation}
    \frac{ \lambda_2}{\lambda_1} \leq 1 +\frac{4}{n}
\end{equation}
and conjectured that
\begin{equation}
    \frac{\lambda_2}{\lambda_1} \leq \left.\frac{\lambda_2}{\lambda_1} \right\vert_{\text{$n$-dimensional ball}} = \left(\frac{j_{n/2}}{j_{n/2-1}}\right)^2,
\end{equation}
where $j_{m}$ is the first positive zero of the Bessel function of order $m$. The original PPW conjecture, along with its $n$-dimensional generalization by Thompson, was proven in 1992 by Ashbaugh and Benguria in \cite{Ashbaugh1991ProofOT, Ashbaugh1992}, which led to a natural question: loosely stated, do more regular shapes maximize the ratio $\lambda_2 / \lambda_1$? In particular, as stated in \cite{numeric}, the polygonal Ashbaugh-Benguria-PPW Conjecture states that the regular $n$-gon in the plane maximizes $\lambda_2 / \lambda_1$ in the class of $n$-gons. More background on the PPW inequality can be found in \cite{henrot2006}.

The purpose of this paper is to show that the ratio $\lambda_2/ \lambda_1$ of eigenvalues of the equilateral triangle is maximized among triangles, as stated below:

\begin{theorem}\label{primary}
For an arbitrary triangle, the following inequality holds:
\begin{equation}
\left.\frac{\lambda_2}{\lambda_1} \right\vert_{\mathrm{triangle}} \leq \left.\frac{\lambda_2}{\lambda_1} \right\vert_{\mathrm{equilateral}} = \frac{7}{3}.
\end{equation}
\end{theorem}

This corresponds to the case $k = 3$ of Conjecture 6.31 in \cite{henrot2017shape} and Conjecture 13 of \cite{numeric}, that is, the triangular case of the polygonal Ashbaugh--Benguria--PPW inequality.

\section{Proof Outline} \label{proof_outline}

\begin{figure}\label{sections}
    \centering
    \includegraphics[scale=0.4]{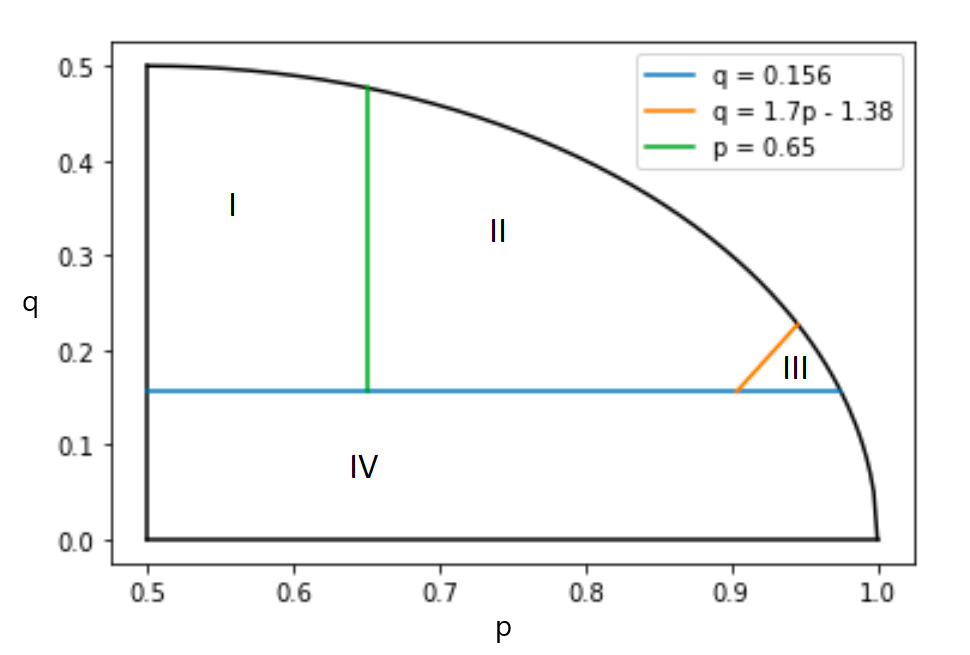}
    \caption{Illustration for single vertex $(p,q)$ of obtuse triangles in cases I, II, III, and IV. The other two vertices are always $(0,0)$ and  $(1,0)$. Note that as $\lambda_2 / \lambda_1$ is invariant under scaling and rigid  motions, we can restrict to this quarter semi-circle without loss of generality.}
    \label{regions}
\end{figure}

In our paper, we prove Theorem \ref{primary} 
by splitting the problem into several cases. 
Our proof of Theorem \ref{primary} relies heavily on work done by Siudeja, who proved in \cite{Siudeja} that Theorem \ref{primary} holds when restricted to acute triangles. Since the acute case is proven in \cite{Siudeja}, we restrict our attention to obtuse and right triangles. We additionally utilize bounds proved by Siudeja and Freitas in \cite{quad}. Once we restrict ourselves to the obtuse case and introduce new bounds for the eigenvalues, we are able to finish the proof of Theorem \ref{primary} with only four  cases, illustrated in Figure \ref{regions}, using mostly simple univariate optimization problems and other elementary techniques. 

In addition for our proof, we build new variational bounds on $\lambda_2$ from those provided in \cite{Siudeja} that are tighter for moderately obtuse triangles, and we apply a simple monotonicity argument to obtain a bound that is effective for very obtuse triangles. We describe this in detail in Section \ref{l2bounds}.



We will use $d$ to denote the diameter of the triangle, which we normalize to 1. We consider triangles in the Euclidean plane with vertices at $(0,0)$, $(1,0)$, and $(p,q)$ without loss of generality. To be right-angled or obtuse at $(p, q)$, the third vertex $(p, q)$ must belong on the boundary of or inside the circle $(p - 1/2)^2 + q^2 = 1/4$. By symmetry, we can focus without loss of generality on the top right quadrant of this circle i.e. when $p \geq 1/2$ and $q \geq 0$. This region is shown in Figure \ref{regions}.

Hence, our triangles have shortest height $h$ equal to $q$ and area $A$ equal to $\frac{q}{2}$. Moreover, we use $\theta$ to denote the smallest angle of the triangle which will be at $(0, 0)$. Thus $\theta = \tan^{-1}(q/p)$.

Our primary strategy is to combine the following estimates for $\lambda_1$ from \cite{quad} and \cite{Siudeja}:
\begin{equation}\label{2.6}
    \lambda_1 \geq \pi^2(1/d +1/h)^2
\end{equation}
and
\begin{equation}\label{2.8}
    \lambda_1 \geq \frac{\theta j_{\pi/\theta}^2}{2A},
\end{equation}
in combination with new bounds on $\lambda_2$. 

We obtain new bounds on $\lambda_2$ using a variational approach with test functions based on known eigenfunctions for the 45-45-90 and 30-60-90 triangles. For very flat triangles, we enclose a rectangle within the triangle. We will refer to these bounds as ``45-45-90 $\lambda_2$ bound," ``30-60-90 $\lambda_2$ bound," and ``rectangle $\lambda_2$ bound."


As seen in Figure \ref{regions}, we divide this region into four areas which we address individually. Area I employs the 45-45-90 bound and bound \eqref{2.6}, Area II uses the 30-60-90 bound and bound \eqref{2.6}, Area III uses the rectangle bound and bound \eqref{2.8}, and finally Area IV employs the rectangle bound and bound \eqref{2.6}.

We now make these $\lambda_2$ estimates precise before going into casework.

\section{Upper Bounds on $\lambda_2$}\label{l2bounds}

Mathematica code reproducing all computations for this section and Section 4 is available on GitHub\footnote{\href{https://github.com/sragavan99/triangle-ppw-inequality}{\texttt{https://github.com/sragavan99/triangle-ppw-inequality}}}.

\subsection{Variational Bounds}



For these bounds on $\lambda_2$ we use the variational characterization

\begin{align*}
    \lambda_2 |_T &= \inf_{f_1, f_2} \sup_\alpha \frac{\int_T |\nabla (\alpha f_1 + f_2)|^2}{\int_T (\alpha f_1 + f_2)^2} \\
    &= \inf_{f_1, f_2} \sup_\alpha \frac{A \alpha ^2 + 2B\alpha + C}{D\alpha^2 + 2E\alpha + F},
\end{align*}
where
\begin{align} 
\begin{aligned} \label{main_integral_coeff}
    A &= \int_T |\nabla f_1|^2, & B &= \int_T \nabla f_1 \cdot \nabla f_2, & C &= \int_T |\nabla f_2|^2, \\
    D &= \int_T f_1^2,  & E &= \int_T f_1f_2,  & F &= \int_T f_2^2.
\end{aligned}
\end{align}
As usual, $f_1, f_2$ must be linearly independent and vanish at the boundary of $T$. To choose test functions $f_1, f_2$, we use the idea of ``transplanting eigenfunctions" used in \cite{quad, hooker_bounds, siudeja_sharp, Siudeja}. We take the first two eigenfunctions of a 45-45-90 or 30-60-90 triangle and transplant them onto $T$ with a suitable affine transformation. These bounds can also be found in \cite{Siudeja}, but the affine transformations used there significantly distort the triangle when it is right or obtuse. We thus obtain better bounds for the obtuse and right cases by choosing different affine transformations that have smaller distortion for right/obtuse triangles; we will point out these differences.
\subsubsection{30-60-90 Bound}

We take our 30-60-90 triangle to have vertices at $(0, 0)$, $(1/2, 0)$, and $(1/2, \sqrt{3}/2)$. On this triangle the first two eigenfunctions are as follows \cite{mccartin}, where for convenience we let $z = \frac{\pi}{3}(2x - 1)$ and $t = \pi (1 - \frac{2y}{\sqrt{3}})$:
\begin{align*}
    \phi_{30, 1}(x, y) &= \sin(4z) \sin(2t) - \sin(5z) \sin(t) - \sin(z) \sin(3t), \\
    \phi_{30, 2}(x, y) &= \sin(5z) \sin(3t) - \sin(2z) \sin(4t) - \sin(7z) \sin(t).
\end{align*}
Let $L_{30}$ be the affine mapping sending $(0, 0)$ to $(1/2, \sqrt{3}/2)$, $(p, q)$ to $(1/2, 0)$, and $(1, 0)$ to $(0, 0)$. This transformation sends the right/obtuse angle of our triangle to the right angle of the 30-60-90 triangle, and it sends the shortest side of the right/obtuse triangle to the shortest side of the 30-60-90 triangle. Thus this preserves the geometry of the triangle reasonably well.

On the other hand, the argument in \cite{Siudeja} starts with the 30-60-90 triangle with vertices at $(0, 0), (1, 0), (0, \sqrt{3})$ and considers an affine mapping preserving $(0, 0)$ and $(1, 0)$ and sending $(p, q)$ to $(0, \sqrt{3})$. For right/obtuse triangles, this is very distortive since the right/obtuse angle at $(p, q)$ is mapped into the $30^{\circ}$ angle at $(0, \sqrt{3})$. Hence we expect our chosen affine mapping to be more effective for the triangles in question.

We take our test functions to be $\phi_{30, 1} \circ L_{30}$ and $\phi_{30, 2} \circ L_{30}$. We can then evaluate coefficients given by \eqref{main_integral_coeff} using these test functions:

\begin{align}
    \begin{aligned}
    A_{30}(p, q) &= \frac{-1594323 + 604800 \pi^2 + 
   4 p (1245184 - 713743 p + 100800 (-3 + 2 p) \pi^2)}{345600 q} \\
   &\qquad + \frac{-2854972 q^2 + 806400 \pi^2 q^2}{345600 q}, \\
   B_{30}(p, q) &= -\frac{2657205 + 4 p (-1507328 + 621593 p) + 2486372 q^2}{354816 q}, \\
   C_{30}(p, q) &= \frac{-1594323 + p (6209536 - 6879600 \pi^2) + 
   28 p^2 (-145849 + 163800 \pi^2)}{1058400 q} \\
   &\qquad + \frac{-4083772 q^2 + 1146600 \pi^2 (3 + 4 q^2)}{1058400 q}, \\
   D_{30}(p, q) &= \frac{3 q} {8}, \quad  E_{30}(p, q) = 0, \quad F_{30}(p, q) = \frac{3 q}{8}.
  \end{aligned}
\end{align}

As expected, $E_{30}(p, q) = 0$ since $\phi_{30, 1}$ and $\phi_{30, 2}$ are orthogonal on the original 30-60-90 triangle and this will be preserved by an affine transformation. Thus our final bound for $\lambda_2$ is

\begin{equation}\label{30bound}
    \lambda_2 \leq \sup_{\alpha} \frac{A_{30}(p, q) \alpha^2 + 2B_{30}(p, q) \alpha + C_{30}(p, q)}{D_{30}(p, q) \alpha^2 + F_{30}(p, q)}.
\end{equation}

\subsubsection{45-45-90 Bound}

We take our 45-45-90 triangle to be that with vertices at $(0, 0), (1, 0), (0, 1)$, for which the first two eigenfunctions are given by:
\begin{align*}
    \phi_{45, 1}(x, y) &= \sin(2\pi x) \sin(\pi y) + \sin(\pi x) \sin(2\pi y), \\
    \phi_{45, 2}(x, y) &= \sin(3\pi x) \sin(\pi y) - \sin(\pi x) \sin(3\pi y).
\end{align*}
These can be derived by noting that eigenfunctions of this triangle can be reflected over the line $y = 1-x$ to obtain an eigenfunction of the unit square that vanishes along this diagonal. We define $L_{45}$ to be the affine mapping sending $(p, q)$ to $(0, 0)$, $(0, 0)$ to $(1, 0)$, and $(1, 0)$ to $(0, 1)$. Note once again that our affine mapping sends the right/obtuse angle at $(p, q)$ to the right angle at $(0, 0)$. In contrast, the work in \cite{Siudeja} (considering the same triangle) works with the affine mapping preserving $(0, 0)$ and $(1, 0)$ and sending $(p, q)$ to $(0, 1)$. Once again, this is very distortive for right/obtuse angles, since the right/obtuse angle at $(p, q)$ is sent to the $45^\circ$ angle at $(0, 1)$. Thus we can expect our affine mapping to yield tighter bounds here as well.

We take our test functions to be $\phi_{45, 1} \circ L_{45}$ and $\phi_{45, 2} \circ L_{45}$. From these test functions, we obtain the following coefficients:

\begin{align}
\begin{aligned}
    A_{45}(p, q) &= \frac{p (256 - 90 \pi^2) + p^2 (-256 + 90 \pi^2) - 256 q^2 + 
 45 \pi^2 (1 + 2 q^2)}{72q}, \\
    B_{45}(p, q) &= \frac{512(1 - 2p)}{175q}, \\
    C_{45}(p, q) &= \frac{5\pi^2(1 - 2p + 2p^2 + 2q^2)}{4q}, \\
 D_{45}(p, q) &= \frac{q}{4}, \quad E_{45}(p, q) = 0, \quad F_{45}(p, q) = \frac{q}{4}. \\
\end{aligned}
\end{align}

Once again, it can be seen without doing any integration that $E_{45}(p, q) = 0$. This gives us the following bound:

\begin{equation}\label{45bound}
    \lambda_2 \leq \sup_{\alpha} \frac{A_{45}(p, q) \alpha^2 + 2B_{45}(p, q) \alpha + C_{45}(p, q)}{D_{45}(p, q) \alpha^2 + F_{45}(p, q)}.
\end{equation}

\subsection{Rectangle Bound}

\begin{figure}
    \centering
    \includegraphics[scale=10.0]{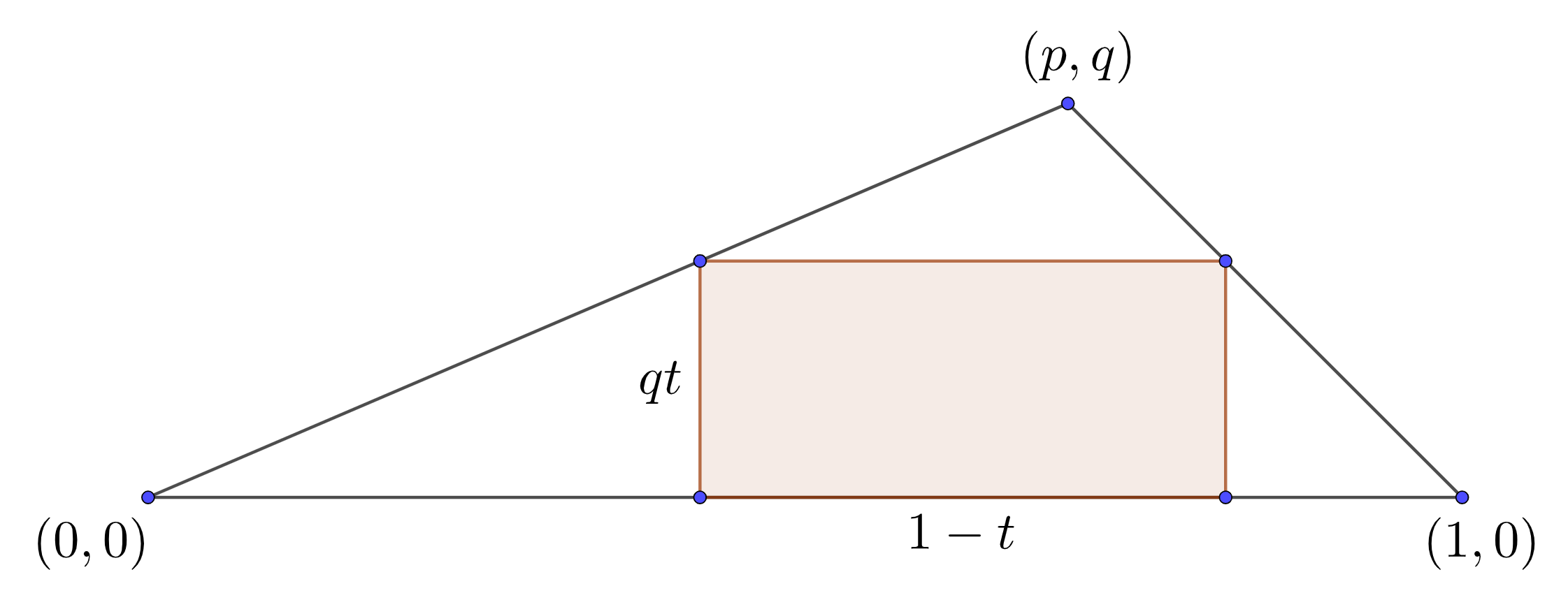}
    \caption{For very obtuse triangles, we bound $\lambda_2$ by enclosing a rectangle within the triangle and using monotonicity.}
    \label{rectangle}
\end{figure}

When our triangle is very obtuse (i.e. $q$ is very small), the bounds on $\lambda_2$ described so far are insufficient. This is not surprising, since in this region our affine transformations are still quite distorted. Thus, we address this case with a different approach. As stated in Section \ref{proof_outline}, this bound is obtained by enclosing a rectangle inside the triangle, with one side aligned with the triangle's diameter. A visualization is given in Figure \ref{rectangle}. As the triangle becomes more obtuse, it becomes closer to the enclosed rectangle in shape, so we expect this estimate to be more effective. It is straightforward to see that if such a rectangle $R$ has height $qt$ for $t \in (0, 1)$, then it will have width $1-t$.

Let us take $t = \frac{1}{1 + \sqrt[3]{4q^2}}$. This clearly is in $(0, 1)$. Moreover, the following inequality holds:

\begin{equation}
    q < 4 \Leftrightarrow q < \sqrt[3]{4q^2} \Leftrightarrow \frac{1}{1+q} > \frac{1}{1+ \sqrt[3]{4q^2}} \Leftrightarrow 1 > (1+q)t \Leftrightarrow 1-t > qt.
\end{equation}

Then by monotonicity of Dirichlet eigenvalues, we obtain the following inequality:

\begin{align*}
    \lambda_2 &\leq \lambda_2(R) = \pi^2(\frac{4}{(1-t)^2} + \frac{1}{(qt)^2}) = \pi^2 \frac{(1 + \sqrt[3]{4q^2})^3}{q^2}.
\end{align*}


This then yields the rectangle bound:

\begin{equation}\label{rectanglebound}
    \lambda_2 \leq \pi^2 \frac{(1 + \sqrt[3]{4q^2})^3}{q^2}.
\end{equation}

Note that a similar bound via rectangle inscription can also be found in \cite{Freitas}.

\section{Proofs in Each Area}

We will now begin to prove Thereom \ref{primary} by splitting into the four cases indicated by Figure \ref{sections}.
\subsection{Area I}

In this region, we have $q \geq 0.156$ and  $p \leq 0.65$, and we employ bounds \eqref{45bound} and \eqref{2.6}. We wish to show that

\begin{equation}
    \frac{A_{45}(p, q) \alpha^2 + 2B_{45}(p, q) \alpha + C_{45}(p, q)}{D_{45}(p, q) \alpha^2 + F_{45}(p, q)} \leq \frac{7}{3} \pi^2 \left(1 + \frac{1}{q}\right)^2
\end{equation}

\noindent for all real $\alpha$ (hence the bound holds for $\sup_\alpha$). Clearing denominators and rearranging, we can equivalently show that the following inequality holds:

\begin{align}\label{45obj}
    \begin{split}
    256 \alpha (288 - 576 p - 175 ((-1 + p) p + q^2) \alpha) - 
 7350 \pi^2 (1 + q)^2 (1 + \alpha^2) \\
 + 7875 \pi^2 (1 + 2 (-1 + p) p + 2 q^2) (2 + \alpha^2) \leq 0.
   \end{split}
\end{align}

For fixed $q$ and  $\alpha$, this is a quadratic in $p$ with leading coefficient given by $31500 \pi^2 - 44800 \alpha^2 + 15750 \pi^2 \alpha^2$. We know that

\begin{equation*}
    31500 \pi^2 - 44800 \alpha^2 + 15750 \pi^2 \alpha^2 > 0
\end{equation*}

\noindent since $15750 \pi^2 > 44800$. For fixed $p$ and  $\alpha$, this is a quadratic in $q$ with leading coefficient $350 (69 \pi^2 - 128 \alpha^2 + 24 \pi^2 \alpha^2)$. We also know that

\begin{equation*}
    350 (69 \pi^2 - 128 \alpha^2 + 24 \pi^2 \alpha^2) > 0
\end{equation*}

\noindent since $24\pi^2 > 128$. Since we wish to prove an upper bound on the LHS of (\ref{45obj}), it suffices to show this upper bound at points where both $p$ and $q$ are extremal assuming the other one is fixed. Thus, we only need to check at the points $(0.5, 0.156), (0.65, 0.156)$, and the points $(p, \sqrt{1/4 - (p - 1/2)^2})$ where $p \in [0.5, 0.65]$.

First, when $(p, q) = (0.5, 0.156)$ we need to show that the following inequality holds:

\begin{equation*}
    \frac{7 (1805312 \alpha^2 - 3 \pi^2 (70267 + 327457 \alpha^2))}{1250} \leq 0.
\end{equation*}

\noindent This clearly holds since $1805312 < 3\pi^2 \cdot 327457$. Second, when $(p, q) = (0.65, 0.156)$, we need to show that

\begin{equation*}
    \frac{128 \alpha (-864000 + 355537 \alpha) - 21 \pi^2 (112318 + 1225453 \alpha^2)}{5000} \leq 0.
\end{equation*}

The LHS is a quadratic in $\alpha$ that attains a maximum of approximately $-1722.58 < 0$ at $\alpha 
\approx -0.265233$, so this inequality holds.

Finally, we deal with the arc at the top of Area I. For convenience, we parametrize the arc as $(1/2 + \sqrt{1/4 - q^2}, q)$ for $q \in [\sqrt{91}/20, 1/2]$. We wish to show that

\begin{equation} \label{semicirc45}
    -73728 \sqrt{1 - 4 q^2} \alpha - 
 525 \pi^2 (-16 - \alpha^2 + 14 q (2 + q) (1 + \alpha^2)) \leq 0.
\end{equation}

For fixed $q$ this is a quadratic in $\alpha$, so it suffices to show that the leading coefficient and the discriminant are both negative. The leading coefficient of \eqref{semicirc45} is given by

\begin{equation*}
    -525 \pi^2 (-1 + 28q + 14q^2),
\end{equation*}

\noindent which on the interval $[\frac{\sqrt{91}}{20}, 0.5]$ achieves a maximum of $\approx -80521.9 < 0$ at $q = \frac{\sqrt{91}}{20}$. Hence, the leading coefficient is negative. The discriminant of \eqref{semicirc45} is given by:

\begin{equation}\label{disc}
    72 (-75497472 (-1 + 4 q^2) - 
   30625 \pi^4 (-8 + 7 q (2 + q)) (-1 + 14 q (2 + q))).
\end{equation}

This is a quartic in $q$, so it may be maximized explicitly; however, we provide a simpler argument here. We claim that it is decreasing in $q$ over the interval $[\sqrt{91}/20, 1/2]$. Indeed its derivative with respect to $q$ is given by:

\begin{equation*}
    -864360000 \pi^4 q^3 - 2593080000 \pi^4 q^2 - 72 (603979776 + 16721250 \pi^4) q + 524790000 \pi^4.
\end{equation*}

This is clearly decreasing for $q >0$, and is hence at most its value at $q = \sqrt{91}/20$ which is $\approx -9.21588 \cdot 10^{10} < 0$.

The discriminant \eqref{disc} is indeed decreasing, and it attains its maximum over this interval at $q = \sqrt{91}/20$ with value $\approx -4.12237 \cdot 10^8 < 0$. This implies that \eqref{semicirc45} holds, completing our proof for Area I.





\subsection{Area II}

In this area, we have $q \geq \max(0.156, 1.7p - 1.38)$ and we utilize equations \eqref{30bound} and \eqref{2.6}. That is, we wish to show that

\begin{equation}
    \frac{A_{30}(p, q) \alpha^2 + 2B_{30}(p, q) \alpha + C_{30}(p, q)}{D_{30}(p, q) \alpha^2 + F_{30}(p, q)} \leq \frac{7}{3} \pi^2\left (1+\frac{1}{q}\right )^2
\end{equation}

\noindent for all real $\alpha$. Upon clearing denominators and rearranging, the problem is to show that the following inequality holds:

\begin{align}\label{big30}
    \begin{split}
        &3((-256 p (-10486784 + 2546775 \pi^2) + 
  28 p^2 (-54958211 + 15523200 \pi^2)\\ &\quad - 
  539 (1594323 + 2854972 q^2) + 54331200 \pi^2 (3 + q (-6 + 5 q))) \alpha^2 \\ &\quad
  +(-2790065250 + 6330777600 p - 2610690600 p^2 - 2610690600 q^2) \alpha\\ &\quad
  -280600848 - 256 p (-4269056 + 4729725 \pi^2) + 
 28 p^2 (-25669424 + 28828800 \pi^2)\\ &\quad - 718743872 q^2 + 7761600 \pi^2 (57 - 42 q + 83 q^2)) \leq 0.
    \end{split}
\end{align}
    
This proof proceeds similarly to our proof in Area I. For fixed $q$ and $\alpha$, \eqref{big30} is a quadratic in $p$ with leading coefficient given by

\begin{equation}\label{30p}
    84 (-25669424 - 7 \alpha (13319850 + 7851173 \alpha) + 2217600 \pi^2 (13 + 7 \alpha^2)),
\end{equation}

and at fixed $p$ and $\alpha$, \eqref{big30} is a quadratic in $q$ with leading coefficient given by

\begin{equation}\label{30q}
    3 (-718743872 - 2610690600 \alpha - 1538829908 \alpha^2 + 7761600 \pi^2 (83 + 35 \alpha^2)).
\end{equation}

As in the proof for Area I, we now show that these leading coefficients are positive for all real $\alpha$. We again check that the discriminant of each is negative, guaranteeing no real roots, which in combination with having a positive leading coefficient implies each is positive for all real $\alpha$. We can bound \eqref{30p} in the following way:
\begin{align*}
    &84(-7 \cdot 7851173 + 7 \cdot 2217600 \pi^2)  > 84(-10^7 + 10^8) > 0.
\end{align*}
\noindent Similarly, we can bound \eqref{30q} in the following way:
\begin{align*}
    3(-1538829908 + 35 \cdot 7761600 \pi^2) > 3(-10^9 + 2 \cdot 10^9) > 0
\end{align*}
\noindent The discriminant corresponding to \eqref{30p} is given by
\begin{align*}
    &(-7 \cdot 84 \cdot 13319850)^2 \\
    &\quad 
    -4 \cdot 84^2 \cdot (-7 \cdot 7851173 + 7 \cdot 2217600\pi^2) \cdot (-25669424 + 13 \cdot 2217600\pi^2) \\
    &\approx -2.4 \cdot 10^{19} < 0,
\end{align*}

\noindent and the discriminant corresponding to \eqref{30q} is given by
\begin{align*}
    &(3 \cdot -2610690600)^2 \\
    &\quad -4\cdot 3^2 \cdot (-718743872 + 7761600 \cdot 83\pi^2) \cdot (-1538829908 + 7761600 \cdot 35\pi^2) \\
    &\approx -7.6 \cdot 10^{20} < 0.
\end{align*}

Therefore, both discriminants of the above polynomials with positive leading coefficients are negative, and they are therefore both always positive for all real $\alpha$. Hence, as in Area I, it suffices to show \eqref{big30} at points where both $p$ and $q$ are extremal assuming the other one is fixed. Thus we only need to check at the point $(0.65, 0.156)$, the line segment covering points $(p, 1.7p - 1.38)$ as $p$ ranges over $[\frac{384}{425}, \frac{1423 + 10\sqrt{1729}}{1945}]$, and the semicircular arc covering points $(p, \sqrt{1/4 - (p - 1/2)^2})$ as $p$ ranges over $[0.65, \frac{1423 + 10\sqrt{1729}}{1945}]$.
\vspace{4mm}

We first address the point $(0.65, 0.156)$. Here we wish to show that:
\begin{align*}
    &\frac{3 (6788089600688 + 7 \alpha (1414193259450 + 1768358569901 \alpha))}{62500} \\
   &\quad + \frac{3 (-7761600 \pi^2 (312007 + 977515 \alpha^2))}{62500} \leq 0. \\
\end{align*}
The LHS is a quadratic in $\alpha$ with leading coefficient $\approx -3.00014 \cdot 10^9 < 0$ and discriminant $\approx -9.6317 \cdot 10^{18} < 0$ so indeed it is always negative.


Next we address the line segment $q = 1.7p - 1.38$. Plugging this into \eqref{big30}, we wish to show that:

\begin{align}\label{lineobj}
\begin{split}
    &\frac{1}{625} (528 (-5857161498 + 15856621390 p - 9928670675 p^2 + 
     11025 (682563\\
     &\quad + 5p (-2211921178 + 1208998385 p)) \alpha + 
  1617 (-4394582298 + 11485165390 p\\
  &\quad - 6941150675 p^2 + 
     126000 (10401 + 5 p (-4566 + 2245 p)) \pi^2) \alpha^2) \leq 0.
\end{split}
\end{align}

For fixed $p$, this is a quadratic in $\alpha$, so we only need to check that its leading coefficient and discriminant are both negative. Firstly, its leading coefficient is:

\begin{align*}
    &\frac{1617}{625} (-4394582298 + 11485165390 p - 6941150675 p^2 \\
   &\quad + 126000 (10401 + 5 p (-4566 + 2245 p)) \pi^2).
\end{align*}

This is a quadratic in $p$, and over the interval of interest it is maximized at $p = \frac{384}{425}$, where its value is $\approx -2.60188 \cdot 10^9 < 0$, so indeed the leading coefficient of \eqref{lineobj} is negative. Next, the discriminant of \eqref{lineobj} is given by:

\begin{align}\label{lineobjdisc}
\begin{split}
&\frac{1764}{390625} (5625 (4620157398 + 5 p (-2211921178 + 1208998385 p))^2 \\
  &\quad - 1936 (-4394582298 + 11485165390 p
  - 6941150675 p^2 + 
     126000 (10401 \\
     &\quad + 5 p (-4566 + 2245 p)) \pi^2) (-5857161498 + 
     15856621390 p
     - 9928670675 p^2 \\
     &\quad + 11025 (682563 + 5 p (-308418 + 171935 p)) \pi^2)).
\end{split}
\end{align}

We now show that \eqref{lineobjdisc} is negative over our interval $[\frac{384}{425}, \frac{1423+10 \sqrt{1729}}{1945}]$. Again, we could do this explicitly but we provide a simpler proof. The second derivative of \eqref{lineobjdisc} is a quadratic in $p$ that is minimized at $\frac{384}{425}$, achieving a minimum of $1.21487 \cdot 10^{21} > 0$, so this quartic is convex over this interval. Hence to show that it is negative it suffices to check that it is negative at the endpoints of our interval. At $p = \frac{384}{425}$, this discriminant is 

\begin{align*}
    &-\frac{1}{20390869140625} ( 63504 \cdot (-1238597349932730480637535561 \\
    &\quad + 1524600 \pi^2 (-147282555087281544521 + %
      24336702382640067000 \pi^2))) \\
      &\approx -4.96593 \cdot 10^{17} < 0,
\end{align*}

while at $p = \frac{1423+10 \sqrt{1729}}{1945}$ it evaluates to 

\begin{align*}
    &-\frac{1}{572451126025} \cdot 142884 \cdot (-93894331197981557997 \cdot (-214373 + 880 \sqrt{1729})\\
   &\quad  + 1355200 \pi^2 \cdot (-114850305 \cdot (-57253311661 + 1065417440 \sqrt{1729})\\
      &\quad  + 21952 \cdot (-155377895789549 + 3421629255040 \sqrt{1729}) \cdot \pi^2))\\
      &\approx -3.44375 \cdot 10^{17} < 0.
\end{align*}


Thus, \eqref{lineobjdisc} is indeed negative. This proves \eqref{lineobj}, achieving the desired result along the line segment $q = 1.7p - 1.38$.


Finally, we show \eqref{big30} along the semi-circular boundary arc; in this domain, we are restricted to $q = \sqrt{1/4 - (p -1/2)^2}$ and $p \in [0.65, \frac{1423+10 \sqrt{1729}}{1945}]$. For convenience, let $p_1 = 0.65$ and $p_2 = \frac{1423+10 \sqrt{1729}}{1945}$. Plugging in $q = \sqrt{1/4 - (p -1/2)^2}$ to ($\ref{big30}$), we want to show that:
\begin{align}\label{semicirc30}
\begin{split}
&27 (-177147 (22 + 7 \alpha) (8 + 77 \alpha) + 236196 p (22 + 7 \alpha) (8 + 77 \alpha) \\
   &\quad + 18110400 p^2 \pi^2 (1 + \alpha^2)
   - 862400 p \pi^2 (73 + 49 \alpha^2) - 2587200 \pi^2 (-19 \\
      &\quad + 14 \sqrt{p(1-p)} + 7 (-1 + 2 \sqrt{p(1-p)}) \alpha^2) \leq 0.
\end{split}
\end{align}
Once again, for fixed $p$ this is a quadratic in $\alpha$. Its leading coefficient is
\begin{align*}
    &27(-95482233 + 127309644 p - 42257600 p \pi^2 + 
   18110400 p^2 \pi^2 \\
   &\quad - 18110400 (-1 + 2 \sqrt{(1 - p) p}) \pi^2) \\ 
   &\coloneqq 27f(p),
\end{align*}
and its constant coefficient is
\begin{align*}
    &27(-31177872 + 41570496 p - 62955200 p \pi^2 +18110400 p^2 \pi^2\\
    &\quad - 2587200(-19 + 14\sqrt{-(-1 + p) p})\pi^2) \\ 
    &\coloneqq 27g(p).
\end{align*}
We show that $f(p)$ is negative in the interval $[p_1, p_2]$ by showing that $f(p)$ is strictly increasing in $[p_1, p_2]$ and that $f(p_2) < 0$. Note that the derivative $f'(p)$ is given by the following:
\[ f'(p) = (127309644 - 42257600\pi^2) + 2\cdot 18110400\pi^2 \cdot p + 18110400 \pi^2\frac{2p-1}{\sqrt{(1-p)p}}. \]
Since $[p_1, p_2] \subset (0.5, 1)$, it trivially holds that both $2p-1$ and $\frac{1}{\sqrt{(1-p)p}}$ are increasing on $[p_1, p_2]$. Since both are also positive on $[p_1, p_2]$, it follows that $\frac{2p-1}{\sqrt{(1-p)p}}$ is increasing on $[p_1, p_2]$. We can then conclude that $f'(p)$ is increasing on $[p_1, p_2]$, and thus:
\begin{align*}
    f'(p) &\ge f'(p_1) \\
    &= (127309644 - 42257600\pi^2) + 2\cdot 18110400\pi^2 \cdot 0.65 + 18110400 \pi^2\frac{0.3}{\sqrt{0.2265}} \\
    &\approx 5.5 \cdot 10^7 > 0.
\end{align*}
Hence $f$ is increasing, and evaluating $f(p_2)$ gives 

\begin{align*}
\begin{split}
    &-95482233 + \frac{127309644 (1423 + 10 \sqrt{1729})}{1945} -  \frac{8451520}{389}(1423 + 10 \sqrt{1729}) \pi^2\\
    &\quad + \frac{724416 (1423 + 10 \sqrt{1729})^2 \pi^2}{151321}\\
    &\quad - 18110400 (-1 + 2 \sqrt{\frac{(1423 + 10 \sqrt{1729}) (1 + \frac{-1423 - 10 \sqrt{1729}}{1945})}{1945}}) \pi^2 \\
    &\approx -1.12135 \cdot 10^8 < 0,
\end{split}
\end{align*}

showing that $f(p)$ is negative for $p\in[p_1,p_2]$, and thus the leading coefficient of \eqref{semicirc30} is negative. It remains to show that its discriminant is also negative.


Similarly, it is easy to see that $g(p)$ is negative for $p\in[p_1,p_2]$. The function $g$ is smooth in the given interval, so its extrema occur at $p_1$, $p_2$, or points within $(p_1,p_2)$ where the derivative of $g$ vanishes. By considering the derivative of $g$ on $(p_1,p_2)$, we find that $g$ has exactly one minimum in $(p_1, p_2)$ and no other local extrema. (Specifically, solving $g' = 0$ in the interval can be reduced to finding the roots of a quartic polynomial, which can be done accurately up to a small error term. Doing so, we find that $g'$ has at most one root in this interval, and we can use the graph of $g'$ and  the intermediate value theorem to see that it has exactly one root.) We call this minimizing value $p_3$ and note that $p_3 \approx 0.81416$. Furthermore $g(p_1) > g(p_3)$ and $g(p_2) > g(p_3)$, so $p_3$ is the unique global minimizer on $[p_1,p_2]$. Thus $g$ is decreasing on $[p_1,p_3]$ and increasing on $[p_3,p_2]$. More details can be found in \texttt{Section4-2.nb} in the GitHub repository.

We have that $-f$ and $-g$ are both decreasing and positive on $[p_3,p_2]$, and hence $fg$ is as well, as is any non-trivial function proportionate to it.

Meanwhile, the square of the linear coefficient, $729 (-310007250 + 413343000 p)^2$, is increasing on $[p_3,p_2]$ and is positive. Thus the discriminant
\[ 729 (-310007250 + 413343000 p)^2 - 2916 f(p) g(p) \] 
is increasing on $[p_3,p_2]$. Since evaluation at $p_2$ gives $\approx -3.4 \cdot 10^{17} < 0$, the discriminant is negative over $[p_3, p_2]$. 

We now check the region $[p_1,p_3]$. In this interval, we have the bounds

\begin{align*}
-2(p- 1/2)^2 +1/2 \le \sqrt{(1-p)p} 
\end{align*}
and
\begin{align*}
\frac{1}{2}-\left(p-\frac{1}{2}\right)^{2}-0.02 \le \sqrt{(1-p)p},
\end{align*}
which hold for $ p \in [p_1, p_3]$. 

Define the following functions, which are modifications of $f$ and $g$ respectively:
\begin{align*}
\begin{split}
    f_0(p) \coloneqq -95482233 + 127309644 p - 42257600 p \pi^2 + 
   18110400 p^2 \pi^2 \\
   - 18110400 (-1 + 2 (\frac{1}{2}-\left(p-\frac{1}{2}\right)^{2}-0.02)) \pi^2,
\end{split}
\end{align*}

\begin{align*}
\begin{split}
    g_0(p) \coloneqq -31177872 + 41570496 p - 62955200 p \pi^2 +18110400 p^2 \pi^2 \\
    - 2587200(-19 + 14(-2(p- 1/2)^2 +1/2))\pi^2.
\end{split}
\end{align*}

Then $f \le f_0$, and additionally $-g$ is positive in $p \in [p_1, p_3]$, so $-fg \le -f_0g $. It is not hard to check that $-f_0$ is positive in $[p_1,p_3]$, and as $g \le g_0$, we similarly have $-f_0 g \le -f_0g_0$. So:
\begin{align*}
729 (-310007250 &+ 413343000 p)^2 - 2916 f(p) g(p) \\
&\le 729 (-310007250 + 413343000 p)^2 - 2916 f_0(p) g_0(p).\\
&=: r(p).
\end{align*}
We will now prove that the LHS is negative throughout $[p_1,p_3]$ by proving this statement for $r(p)$. To do this, we first want to prove as a lemma that $r$ is concave up in $[p_1,p_3]$; $r$ is a quartic polynomial, meaning that its second derivative $r''$ is quadratic. To find the critical point of the quadratic $r''$, we find the root of its derivative. The derivative is given by $-344307200786841600000\pi^4 p -241212414904592947200 \pi^2
    + 253038466610012160000\pi^4,$
which has its root at $p = (-17301357 + 18149600 \pi^2)/(24696000 \pi^2) \approx 0.663938$, achieving the value $\approx 1.32883 \cdot 10^{21} > 0$. Since $r''$ is a quadratic whose leading coefficient is $-172153600393420800000\pi^4 < 0$, $r''(p)$ is a maximum, and since $p \in [p_1,p_3]$, we merely have to check that both $r''(p_1) > 0$ and $r''(p_3) > 0$ to conclude that $r''(p) > 0$ for all $p \in [p_1,p_3]$. It suffices to prove that $r''(q) > 0$ for some $q > p_3$ in place of the statement $r''(p_3) > 0$ to avoid approximations involving the evaluation of $r''(p_3)$. At $p=0.65$ and $p = 0.83 > p_3$, $r''$ evaluates to $\approx 1.32557\cdot10^{21}$ and $\approx 8.66388\cdot10^{20}$, respectively. Therefore, $r''(p) > 0$ for all $p \in [p_1,p_3]$, and $r$ is concave upwards in this region.

Thus it suffices to check that $r(p) < 0$ at $p = p_1$ and $p = 0.83 > p_3$. At these points, $r(p)$ has the values $\approx -4.39\cdot10^{18} < 0$ and $\approx -1.62 \cdot 10^{18} < 0$, respectively. Hence, $r(p) < 0$ for all $p \in [p_1,p_3]$, implying that the same holds for the discriminant of \ref{semicirc30}.

This completes the proof of Equation \eqref{big30}, and hence Theorem \ref{primary} over Area II.

\subsection{Area III}

\begin{figure}
    \centering
    \includegraphics[scale=1.3]{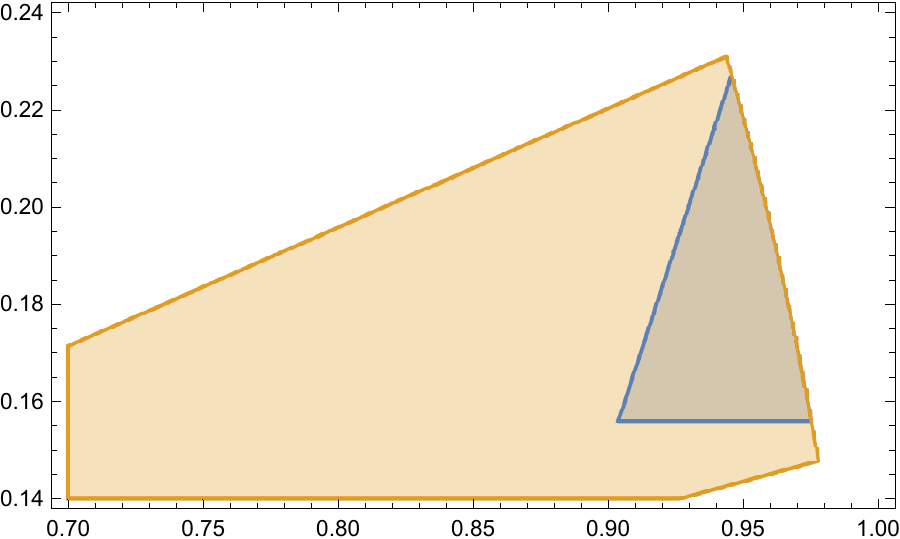}
    \caption{Visualization showing that the range of angles $\theta$ in Area III is contained in $[0.15, 0.24]$. Area III (shown in blue) is covered by the region $p\tan (0.15) \leq q \leq p \tan(0.24)$ (shown in orange).}
    \label{angrange}
\end{figure}

For this area, we use bounds \eqref{rectanglebound} and \eqref{2.8}. Here we are concerned with the region where $0.156 \leq q \leq 1.7p - 1.38$. Note that in this region the range of angles $\theta$ achieved is contained in $[0.15, 0.24]$. This is straightforward to check since all the conditions of interest are linear or quadratic inequalities in $p, q$, and is shown visually in Figure \ref{angrange}. Detailed code for producing this figure and the other computations in this section can be found in \texttt{Section4-3.nb} in the GitHub repository.

Noting that $A = \frac{q}{2}$, we can then rearrange our target inequality to the following:
\begin{align*}
    \frac{\pi^2\frac{(1 + \sqrt[3]{4q^2})^3}{q^2}}{\frac{\theta j_{\pi/\theta}^2}{2A}} \leq \frac{7}{3}
    \Leftrightarrow 3\pi^2\frac{(1 + \sqrt[3]{4q^2})^3}{q^2} &\leq 7\frac{\theta j_{\pi/\theta}^2}{2A} \\
    \Leftrightarrow 3\pi^2\frac{(1 + \sqrt[3]{4q^2})^3}{q} &\leq 7\theta j_{\pi/\theta}^2.
\end{align*}

Call the LHS $f(q)$ and the RHS $g(\theta)$. We claim that these are decreasing functions of $q$ and $\theta$ respectively in the region we are interested in. We first check this for $f(q)$, omitting the constant factor of $3\pi^2$:

\begin{align*}
    \frac{d}{dq} \frac{(1 + \sqrt[3]{4q^2})^3}{q} &= \frac{2 \cdot 2^{2/3}(1 + \sqrt[3]{4q^2})^2}{q^{4/3}} - \frac{(1 + \sqrt[3]{4q^2})^3}{q^2} \\
    &= \frac{(1 + \sqrt[3]{4q^2})^2}{q^2} (2 \cdot 2^{2/3} q^{2/3} - (1 + \sqrt[3]{4q^2})) \\
    &= \frac{(1 + \sqrt[3]{4q^2})^2}{q^2} (\sqrt[3]{4q^2} - 1) \leq 0,
\end{align*}
since $q \leq 1/2$.

Now for $g(\theta)$, let $t = \pi/\theta \in [13, 21]$. Since $t$ is decreasing with respect to $\theta$, we wish to show (omitting the constant factor of 7) that $\frac{j_t^2}{t}$ is increasing with respect to $t$ in this region:
\begin{align*}
    \frac{d}{dt} \frac{j_t^2}{t} &= \frac{t2j_t \frac{d j_t}{dt} - j_t^2}{t^2}.
\end{align*}
So we want to show that $t2j_t \frac{d j_t}{dt} - j_t^2 \geq 0 \Leftrightarrow j_t^2 \leq 2tj_t \frac{d j_t}{dt} \Leftrightarrow j_t \leq 2t \frac{dj_t}{dt}$. To do this, we use the following results shown in \cite{elbert}:
\begin{enumerate}
    \item $\frac{dj_t}{dt} > 1$ (this is Lemma 1.1 of \cite{elbert}, which is applicable here since $t > 0$ and $j_0 \approx 2.40 > \frac{1}{4}$.
    \item $j_t$ is concave as a function of $t$ (this is Corollary 3.3 of \cite{elbert}).
\end{enumerate}
The first point means we just need to show that $j_t \leq 2t$. Note that at $t = 13$ we have $j_t \leq 17.802 < 26 = 2t$ so it suffices to show that $\frac{dj_t}{dt} \leq 2$ for $t \geq 13$. By the second point, we know that $\frac{dj_t}{dt}$ is non-increasing so we just need to show that $\left. \frac{dj_t}{dt}\right\vert_{t = 13} \leq 2$. But this is straightforward; by concavity, the LHS is at most $j_{13} - j_{12} \approx 1.10 < 2$.

\begin{figure}
    \centering
    \includegraphics[scale=1]{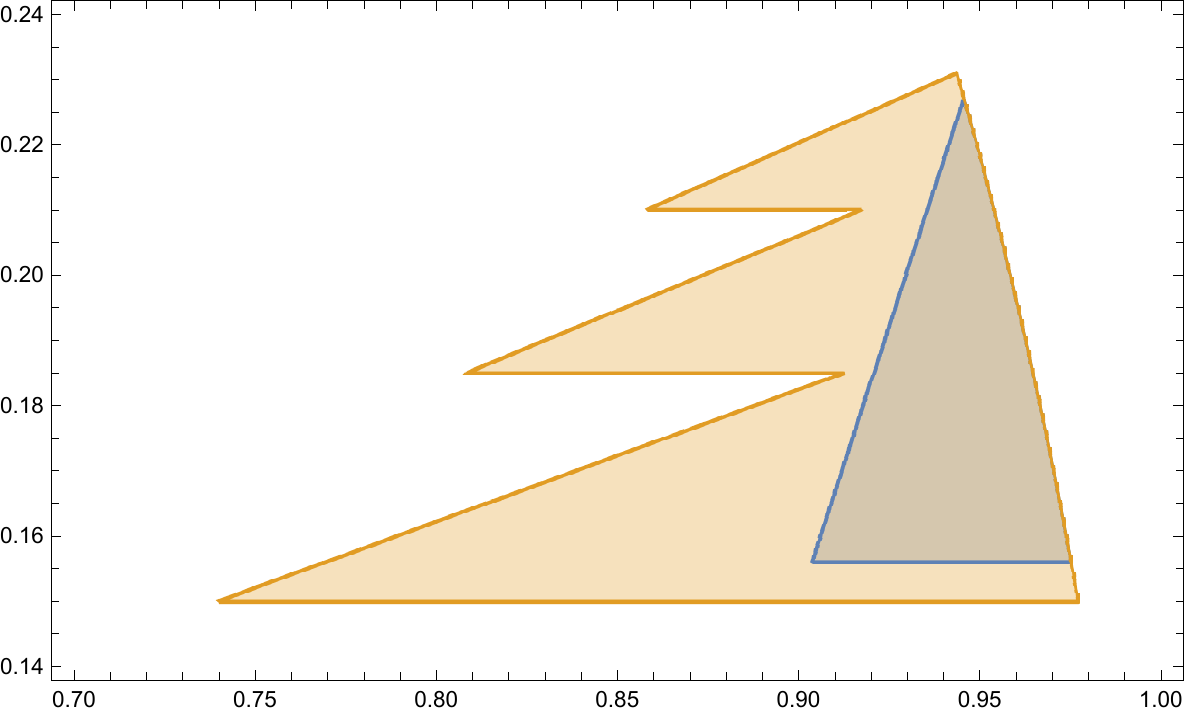}
    \caption{Final step for the proof in Area III (shown in blue here). The orange region is the union of three sector-like regions, each of which corresponds to one $S_{q_0, \theta_0}$. As desired, the orange region covers Area III.}
    \label{regionscase3}
\end{figure}

Hence our claim is proven. Then the key observation is this, if we have a particular pair $(q_0, \theta_0)$ such that $f(q_0) \leq g(\theta_0)$, then whenever $q \geq q_0$ and $\theta \leq \theta_0$ we have $f(q) \leq f(q_0) \leq g(\theta_0) \leq g(\theta)$. So if we call the set of points inside our quarter circle satisfying $q \geq q_0$ and $\theta \leq \theta_0$ as $S_{q_0, \theta_0}$ then it suffices to specify a set of pairs $(q_0, \theta_0)$ each satisfying $f(q_0) \leq g(\theta_0)$ such that the sets $S_{q_0, \theta_0}$ collectively cover Area III. We take the following three pairs:
\begin{itemize}
    \item $q_0, \theta_0 = 0.15, 0.2$, satisfies $f(q_0)/g(\theta_0) \approx 0.9929 < 1$,
    \item $q_0, \theta_0 = 0.185, 0.225$, satisfies $f(q_0)/g(\theta_0) \approx 0.9943 < 1$,
    \item $q_0, \theta_0 = 0.21, 0.24$, satisfies $f(q_0)/g(\theta_0) \approx 0.9959 < 1$.
\end{itemize}
It is straightforward to verify that the three $S_{q_0, \theta_0}$'s thus defined cover the region of interest (since all the conditions of interest are linear or quadratic inequalities in $p, q$). This is visually shown in Figure \ref{regionscase3}.

\subsection{Area IV}

As mentioned earlier, here we use bounds \eqref{rectanglebound} and \eqref{2.6}. We are concerned with the region $q \leq 0.156$, where we wish to show that:

\begin{align*}
    \frac{\pi^2 \frac{(1 + \sqrt[3]{4q^2})^3}{q^2}}{\pi^2(1 + \frac{1}{q})^2} &\leq \frac{7}{3} \\
    \Longleftrightarrow \frac{(1 + \sqrt[3]{4q^2})^3}{(q+1)^2} &\leq \frac{7}{3} \\
    \Longleftrightarrow 3(1 + \sqrt[3]{4q^2})^3 - 7(q+1)^2 &\leq 0.
\end{align*}

Call the LHS $f(q)$. Note that $f(0.156) \approx -0.0177 < 0$ so it suffices to show that $f$ is non-decreasing on $[0, 0.156]$. To do this, substitute $x = \sqrt[3]{q}$ (note that this is an increasing function of $q$) and differentiate $f$ with respect to $x$. We want to show that this derivative is non-negative on $[0, \sqrt[3]{0.156}]$, i.e., that $6x(3 \cdot 2^{2/3} - 7 x + 12 \cdot 2^{1/3} x^2 + 5 x^4) \geq 0$. But this is clear since $x \geq 0$ and $x \leq \sqrt[3]{0.156} < \frac{3 \cdot 2^{2/3}}{7} \Rightarrow 3 \cdot 2^{2/3} - 7x > 0$. This completes our proof for this case.

\section*{Acknowledgements}
We wish to thank Javier G\'omez-Serrano for introducing this problem to us in his class and for guiding us while writing this paper. We also thank the Princeton University Department of Mathematics.

\bibliography{refer}

\providecommand{\bysame}{\leavevmode\hbox to3em{\hrulefill}\thinspace}
\providecommand{\MR}{\relax\ifhmode\unskip\space\fi MR }
\providecommand{\MRhref}[2]{%
  \href{http://www.ams.org/mathscinet-getitem?mr=#1}{#2}
}
\providecommand{\href}[2]{#2}
\begin{thebibliography}{10}

\bibitem{numeric}
P.~Antunes and P.~Freitas, \emph{A numerical study of the spectral gap}, J.
  Phys. A: Math. Theor. \textbf{41} (2008).

\bibitem{Ashbaugh1991ProofOT}
Mark~S. Ashbaugh and Rafael~D. Benguria, \emph{Proof of the
  {P}ayne-{P}{\'o}lya-{W}einberger conjecture}, Bulletin of the American
  Mathematical Society \textbf{25} (1991), 19--29.

\bibitem{Ashbaugh1992}
\bysame, \emph{A sharp bound for the ratio of the first two eigenvalues of
  {D}irichlet {L}aplacians and extensions}, The Annals of Mathematics
  \textbf{135} (1992), no.~3, 601.

\bibitem{elbert}
Arpad Elbert, Luigi Gatteschi, and Andrea Laforgia, \emph{On the concavity of
  zeros of {B}essel functions}, Applicable Analysis \textbf{16} (1983), no.~4,
  261--278 (en).

\bibitem{Freitas}
Pedro Freitas, \emph{Precise bounds and asymptotics for the first dirichlet
  eigenvalue of triangles and rhombi}, Journal of Functional Analysis
  \textbf{251} (2007), 376--398.

\bibitem{quad}
{Freitas, Pedro} and {Siudeja, Bartłomiej}, \emph{Bounds for the first
  {D}irichlet eigenvalue of triangles and quadrilaterals}, ESAIM: COCV
  \textbf{16} (2010), no.~3, 648--676.

\bibitem{henrot2006}
A.~Henrot, \emph{Extremum problems for eigenvalues of elliptic operators},
  Birkhäuser Basel, 2006.

\bibitem{henrot2017shape}
\bysame, \emph{Shape optimization and spectral theory}, De Gruyter Open, 2017.

\bibitem{hooker_bounds}
W.~Hooker and M.~H. Protter, \emph{Bounds for the first eigenvalue of a rhombic
  membrane}, Journal of Mathematics and Physics \textbf{39} (1960), no.~1-4,
  18--34 (en).

\bibitem{mccartin}
Brian~J. McCartin, \emph{Eigenstructure of the {Equilateral} {Triangle}, {Part}
  {I}: {The} {Dirichlet} {Problem}}, SIAM Review \textbf{45} (2003), no.~2,
  267--287 (en).

\bibitem{PPW}
L.~E. Payne, G.~Pólya, and H.~F. Weinberger, \emph{On the ratio of consecutive
  eigenvalues}, Journal of Mathematics and Physics \textbf{35} (1956), no.~1-4,
  289--298.

\bibitem{siudeja_sharp}
Bartłomiej Siudeja, \emph{Sharp bounds for eigenvalues of triangles}, The
  Michigan Mathematical Journal \textbf{55} (2007), no.~2, 243--254 (en).

\bibitem{Siudeja}
Bartłomiej Siudeja, \emph{Isoperimetric inequalities for eigenvalues of
  triangles}, Indiana University Mathematics Journal \textbf{59} (2010), no.~3,
  1097--1120.

\bibitem{thompson}
C.~J. Thompson, \emph{On the ratio of consecutive eigenvalues in
  \textit{N}-dimensions}, Stud. Appl. Math \textbf{48} (1969).

\end{thebibliography}

\end{document}